\documentclass[12pt]{article}
\usepackage{}
\usepackage{amssymb}
\usepackage{mathrsfs}
\usepackage{amsmath,amssymb,theorem}
\usepackage{graphicx}
\usepackage{subfigure}
\usepackage{psfrag}
\usepackage{color}
\usepackage{enumerate}
\newtheorem{thm}{Theorem}[section]

\newtheorem{lem}[thm]{Lemma}
\theorembodyfont{\rmfamily}

\def\less{\backslash}

\def\qed{ \hfill $\square$}

\title{Improved bounds for rainbow numbers of matchings in plane triangulations}

\author{\small {Zhongmei Qin$^1$, Yongxin Lan$^2$\footnote{The corresponding author.}, Yongtang Shi$^2$}\\
{\small $^1$School of Science}\\
{\small Chang'an University, Xi'an, Shaanxi 710064, P.R. China}\\
{\small  $^2$Center for Combinatorics and LPMC}\\
{\small Nankai University, Tianjin 300071, P.R. China}\\
{\small Emails: qinzhongmei90@163.com, yxlan0@126.com, shi@nankai.edu.cn}\\
}
\date{}
\begin{document}\maketitle
\begin{abstract}
Given two graphs $G$ and $H$, the {\it rainbow number} $rb(G,H)$ for $H$ with respect to $G$ is defined as the minimum number $k$ such that any $k$-edge-coloring of $G$ contains a rainbow $H$, i.e., a copy of $H$, all of whose edges have different colors. Denote by $kK_2$ a matching of size $k$ and $\mathcal {T}_n$ the class of all plane triangulations of order $n$, respectively. In [S. Jendrol$'$, I. Schiermeyer and J. Tu, Rainbow numbers for matchings in plane triangulations, Discrete Math. 331(2014), 158--164], the authors determined the exact values of $rb(\mathcal {T}_n, kK_2)$ for $2\leq k \le 4$ and proved that $2n+2k-9 \le rb(\mathcal {T}_n, kK_2) \le 2n+2k-7+2\binom{2k-2}{3}$ for $k \ge 5$. In this paper, we improve the upper bounds and prove that $rb(\mathcal {T}_n, kK_2)\le 2n+6k-16$ for $n \ge 2k$ and $k\ge 5$. Especially, we show that $rb(\mathcal {T}_n, 5K_2)=2n+1$ for $n \ge 11$. \\[2mm]
\end{abstract}

\section{Introduction}

All graphs in this paper are undirected, finite and simple. We follow \cite{BM08} for graph theoretical notation and terminology not defined here. Let $G$ be a connected graph with vertex set $V(G)$ and edge set $E(G)$. For any two disjoint subsets $X$ and $Y$ of $V(G)$, we use $E_G(X,Y)$ to denote the set of edges of $G$ that have one end in $X$ and the other in $Y$. 
We also denote $E_G(X,X)=E_G(X)$.
Let $e(G)$ denote the number of edges of $G$, $e_G(X,Y)$ the number of edges of $E_G(X,Y)$, $e_G(X)$ the number of edges of $E_G(X)$.
If $X=\{x\}$, then we write $E_G(x,Y)$ and $e_G(x,Y)$, respectively. For a vertex $x\in V(G)$, We use $N_G(x)$ to denote the set of vertices in $G$ which are adjacent to $x$. We define $d_G(x)=|N_G(x)|$. Given vertex sets $X, Y \subseteq V(G)$,
the subgraph of $G$ induced by $X$, denoted $G[X]$, is the graph with vertex set $X$ and edge set $\{xy \in E(G) : x, y \in X\}$. We denote by $Y \less X$ the set $Y - X$.

A subgraph of an edge-colored graph is {\it rainbow} if all of its edges are colored distinct. Given two graphs $G$ and $H$, the {\it rainbow number} $rb(G,H)$ for $H$ with respect to $G$ is defined as the minimum number $k$ such that any $k$-edge-coloring of $G$ contains a rainbow copy of $H$. When $G=K_n$, the rainbow number is closely related to anti-Ramsey number, which was introduced by Erd\H{o}s, Simonovits and S\'{o}s \cite{ESS75} in 1975. The {\it anti-Ramsey number}, denoted by $f(K_n,H)$, is the maximum number $c$ for which there is a way to color the edges of $K_n$ with $c$ colors such that every subgraph $H$ of $K_n$ has at least two edges of the same color. Clearly, $rb(K_n, H)=f(K_n, H)+1$.

Let $\mathcal{T}_n$ denote the class of all plane triangulations of order $n$. We denote by $rb(\mathcal{T}_n,H)$ the minimum number of colors $k$ such that, if $H\subseteq T_n\in \mathcal{T}_n$, then any edge-coloring of $T_n$ with at least $k$ colors contains a rainbow copy of $H$. The rainbow number has been widely studied. The rainbow numbers for matchings with respect to complete graph has been completely determined step by step in \cite{CLT09,ESS75,FKSS09, S04}. Also, the rainbow numbers for some other special graph classes in complete graphs have been obtained, see \cite{A83,HY12,J02,JW03,J17,JL09,MN02}. Meanwhile, the researchers studied the rainbow number when host graph changed from the complete graph to others, such as complete bipartite graphs (\cite{AJK04,LTJ09}), planar graphs (\cite{HJSS15,JST14,JYNC,LSS17}), hypergraphs (\cite{OY13}), etc. For more results on rainbow numbers, we refer to the survey \cite{FMO10}.

In this paper we study the rainbow number when host graphs are plane triangulations. Let $\mathcal {T}_n$ be the family of all plane triangulations on $n$ vertices. As one of the most important structures in graphs, the study of rainbow number in plane triangulations $rb(\mathcal {T}_n, H)$ was initiated by Hor\v{n}\'{a}k et al. \cite{HJSS15}. Hor\v{n}\'{a}k et al. \cite{HJSS15} investigated the rainbow numbers for cycles. Very recently, Lan, Shi and Song \cite{LSS17} improve some bounds for the rainbow number of cycles, and also get some results for paths.  Jendrol$'$, Schiermeyer and Tu \cite{JST14} investigated the rainbow numbers for matchings in plane triangulations. We summarize their results as follows, where $kK_2$ denote a matching of size $k$.

\begin{thm}[\cite{JST14}]\label{th2}

\begin{description}
\item[(1)] $2n+2k-9 \le rb(\mathcal {T}_n, kK_2) \le 2n+2k-7+2\binom{2k-2}{3}$ for all $k \ge 5$.
\item[(2)] $rb(\mathcal {T}_n, 2K_2)=\begin{cases}
                         4 & n=4;\\
                         2 & n\ge 5.
                         \end{cases}$
\item[(3)] $rb(\mathcal {T}_n, 3K_2)=\begin{cases}
                         8 & n=6;\\
                         n+1 & n\ge 7.
                         \end{cases}$
\item[(4)] $rb(\mathcal {T}_n, 4K_2)=2n-1$ for all $n \ge 8$.
\end{description}
\end{thm}

Recently, Jin and Ye \cite{JYNC} investigated the rainbow numbers of $kK_2$ in the maximal outerplanar graphs.
In this paper, we improve the upper bounds and prove that $rb(\mathcal {T}_n, kK_2)\le 2n+6k-16$ for $n \ge 2k$ and $k\ge 5$. Especially, we show that $rb(\mathcal {T}_n, 5K_2)=2n+1$ for $n \ge 11$ by using the method of Jendrol$'$, Schiermeyer and Tu \cite{JST14}.

\begin{thm}\label{them2}
For $n \ge 2k$ and $k\ge 5$,
$rb(\mathcal {T}_n, kK_2)\le 2n+6k-16$.
\end{thm}

\begin{thm}\label{them1}
For $n \ge 11$, $rb(\mathcal {T}_n, 5K_2)=2n+1$.
\end{thm}

The following theorem will be used in our proof. A graph $G$ is called {\it factor-critical} if $G-v$ contains a perfect matching for each $v \in V(G)$. A graph is called {\it hypoHamiltonian} if for every vertex $u$, $G-u$ is Hamiltonian.

\begin{thm}[\cite{LP86}]\label{th1}
Given a graph $G=(V,E)$ and $|V|=n$, let $d$ be the size of a maximum matching of $G$. Then there exists a subset $S$ with $|S|\le d$ such that
$$d=\frac{1}{2}(n-(o(G-S)-|S|)),$$
where $o(H)$ is the number of components in the graph $H$ with an odd number of vertices. Moreover, each odd component of $G-S$ is factor-critical.
\end{thm}

\begin{lem}\label{lem}
Let $G$ be a planar triangulation on $n\ge 4$ vertices. Then
\begin{enumerate}[(a)]

  \item  (\cite{J81,JST14}) for $5\le n\le7$, $G$ is hypoHamiltonian.
  \item $G$ is 3-connected.
\end{enumerate}

\end{lem}

\section{Proof of Theorem \ref{them2}.}
By induction on $k$. The statement is true for $k\leq 4$ by Theorem \ref{th2}. Now we assume $k\ge 5$. Let $T_n$ be a plane triangulation on $n$ vertices. By contradiction, let $c$ be an edge-coloring of $T_n$ with $2n+6k-16$ colors such that $T_n$ does not contain any rainbow $kK_2$. Let $G$ be a rainbow spanning subgraph of $T_n$ with $2n+6k-16$ edges. Then $G$ is $kK_2$-free. Since $2n+6k-16>2n+6(k-1)-16$, $G$ contains a $(k-1)K_2$ by the induction hypothesis.
Let $u_1w_1, u_2w_2, \ldots, u_{k-1}w_{k-1}$ be a $(k-1)K_2$ of $G$, and let $H$ be an induced subgraph by $\{u_1,\ldots,u_{k-1}, w_1,\ldots, w_{k-1}\}$ in $G$. Then $e_G(H)\le 3(2k-2)-6=6k-12$.

Let $R=V(G)\less V(H)$. Since $G$ is $kK_2$-free, $E(G[R])=\emptyset$.
Then we have $G-E_G(H)$ is a bipartite planar graph with $n$ vertices, which implies $e_G(V(H),R)\le 2n-4$. Thus, $e(G)= e_G(H)+e_G(V(H),R)\le 6k-12+2n-4=2n+6k-16$. Since $e(G)=2n+6k-16$, we have $e_G(H)=6k-12$ and $e_G(V(H),R)=2n-4$. Hence, $G[V(H)]$ is a plane triangulation with $2k-2$ vertices and $G-E_G(H)$ is a maximal bipartite planar graph with $n$ vertices.
Since $u_1w_1\in E(G)$, there must exist a quadrangular face with vertices $u_1,r_1,w_1,r_2$ in order in $G-E_G(H)$, where $r_1,r_2\in R$. But then the graph induced by the edges $u_1r_1$, $w_1r_2$, $u_2w_2$,\ldots, $u_{k-1}w_{k-1}$ is a rainbow subgraph of $T_n$ isomorphic to $kK_2$, a contradiction.

\section{Proof of Theorem \ref{them1}.}

By Theorem \ref{th2}, we only need to show that $rb(\mathcal {T}_n, 5K_2)\le 2n+1$. Suppose $rb(\mathcal {T}_n, 5K_2)\ge 2n+2$. Then there exists a plane triangulation $T_n$ on $n$ vertices containing no  rainbow $5K_2$ under an edge-coloring $c$ used $2n+1$ colors. Let $G \subset T_n$ be a rainbow spanning subgraph with $2n+1$ edges. Then $G$ has no a copy of $5K_2$. By Theorem \ref{th2}, $G$ has a copy of $4K_2$.
By Theorem \ref{th1}, there exists an $S\subseteq V(G)$ with $|S|=s\le 4$, such that $q=o(G-S)=n-8+s$. Let $A_1, \ldots, A_q$ be all the odd components of $G-S$. Assume $|V(A_i)|=a_i$ for each $i\in[q]$
and $a_1 \ge a_2 \ge \cdots \ge a_q$. Let $t=\min\{i: a_i=1\}$ and $V_0=\{v_t,\ldots,v_q\}$, where $v_j\in V(A_j)$. Assume $d_G(v_t)\ge d_G(v_{t+1})\ge \cdots \ge d_G(v_q)$. Let $B$ denote the set of vertices of all the even components of $G-S$. We first prove a useful claim.\medskip

\noindent {\bf Claim.} If $G$ contains two edge-disjoint $4K_2$, say $M_1$, $M_2$, then $E_{T_n}(V(G)\less V(M_1\cup M_2))=\emptyset$.\medskip

\noindent {\bf Proof.} Suppose $e\in E_{T_n}(V(G)\less V(M_1\cup M_2))$. Assume that $c(e)=c(e')$ for some $e'\in M_1$. Then $e\cup M_2$ is a rainbow $5K_2$, a contradiction.\qed\medskip

Assume first $s\in\{0,1\}$. Let $G'$ be a copy of $G$ with $a_1=(n-s)-q+1=9-2s$ and $a_2=\cdots=a_q=1$.  It is easy to check that $e(G)\le e(G')\le 3(a_1+s)-6+s(q-1)=sn+(s^2-12s+21)<2n+1$ since $n\ge 11$, a contradiction. Hence, $s\in\{2,3,4\}$.\medskip

Suppose $s=2$. Then $q=n-6$. Let $S=\{w_1,w_2\}$. We claim $B=\emptyset$. Suppose $B\ne \emptyset$. Then $|B|\le n-q-s=4$ and so either $|B|=2$ or $|B|=4$. If $|B|=2$, then $a_1=3$ and $a_i=1$ for each $i\ge2$. Hence, $e(G)=e(G[S\cup V(A_1)])+e_G(S,V_0\cup B)+e_G(B)\le (3\cdot 5-6)+2(n-5)+1<2n-1$, a contradiction. Thus $|B|=4$ and so $a_i=1$ for each $i\in[q]$. Then $e(G)=e(G[S\cup B])+e_G(S,V_0)=(3\cdot 6-6)+2(n-6)<2n+1$, a contradiction. Hence, $B=\emptyset$ and so $a_1\in\{3,5\}$. Suppose $a_1=5$. Then $a_i=1$ for each $i\ge 2$. We see $2n+1=e(G)=e(G[S\cup V(A_1)])+e_G(S,V_0)\le (3\cdot 7-6)+2(n-7)=2n+1$, which implies $G[S\cup V(A_1)]$ is a plane triangulation and $d_G(x)=2$ for each $x\in V_0$. Note that $|V_0|\ge 4$ as $n\ge11$. By Lemma~\ref{lem}(a), $G- \{v_i,v_j\}$ for any $v_i,v_j\in V_0$ contains two edge-disjoint $4K_2$. By Claim, $E_{T_n}(V_0)=\emptyset$. Since $\delta(T_n)\ge 3$, $e_{T_n}(x,V(A_1))\ge 1$ for each $x\in V_0$.
But then $T_n$ contains a $K_{3,3}$-minor (with one part $\{v_2,v_3,v_4\}$ and the other part $\{w_1,w_2,V(A_1)\}$), a contradiction. Thus $a_1=3$, and so $a_2=3$ and $a_i=1$ for each $i\ge3$. By Lemma~\ref{lem}(b), $e(G[S\cup V(A_1)\cup V(A_2)])\le 3\cdot 8-7$. Thus, $2n+1=e(G)=e(G[S\cup V(A_1)\cup V(A_2)])+e_G(S,V_0))\le (3\cdot 8-7)+2(n-8)=2n+1$, which implies that $e(G[S\cup V(A_1)\cup V(A_2)])=17$ and $d_G(x)=2$ for each $x\in V_0$.
Let $V(A_1)=\{u_1,u_2,u_3\}$ and $V(A_2)=\{u_4,u_5,u_6\}$. By Theorem~\ref{th1}, $A_1=A_2=K_3$. We may assume $w_1u_1,w_2u_4\in E(G)$. Then $G-\{v_i,v_j\}$ for any $v_i,v_j\in V_0$ contains two edge-disjoint $4K_2$, say $M_1=\{w_1u_1,w_2v_k,u_2u_3,u_4u_5\}$ and $M_2=\{w_2u_4,w_1v_k,u_5u_6,u_1u_2\}$, where $v_k\in V_0\less \{v_i,v_j\}$. By Claim, $E_{T_n}(V_0)=\emptyset$.
Since $\delta(T_n)\ge 3$, $e_{T_n}(x,V(A_1)\cup V(A_2))\ge 1$ for each $x\in V_0$. Thus, $e_{T_n}(V(A_1),V(A_2))=0$ since otherwise $T_n$ contains a $K_{3,3}$-minor (with one part $\{v_3,v_4,v_5\}$ and the other part $\{w_1,w_2,V(A_1)\cup V(A_2)$), and $e_{T_n}(A_i,V_0)\ge 1$ for all $i\in[2]$ since otherwise $T_n$ contains a $K_{3,3}$-minor (with one part $\{v_3,v_4,v_5\}$ and the other part $\{w_1,w_2,V(A_{3-i})\}$).
We claim that $e_G(S,A_i)\ge 5$ for each $i\in[2]$. Assume that $e_G(S,A_1)\le 4$. We see $e_G(S,A_1)=4$, $e_G(S,A_2)=6$  and $w_1w_2\in E(G)$ because $e_G(S,A_2)\le 6$ and $e(G[S\cup V(A_1)\cup V(A_2)])=17$. But then $e(G[S\cup V(A_2)])=10>3\cdot 5-6$, a contradiction.
Moreover, each vertex in $V_0$ has no common neighbors in $V(A_1)\cup V(A_2)$. Suppose there exists two vertices, say $v_3,v_5$ such that $v_3u_5,v_5u_5\in E(T_n)$. Then $T_n$ contains a $K_{3,3}$-minor (with one part $\{v_3,v_5,\{u_4,u_6\}\}$ and the other part $\{w_1,w_2,u_5\}$).  Hence, we further assume $v_3u_5,v_4u_2,v_5u_6\in E(T_n)$. We claim that $c(v_3u_5)=c(u_4u_6)$. Suppose not, since $e_G(S,A_1)\ge 4$, we may assume $\{u_1,u_2\}\subset N_{A_1}(w_1)$. Then $G-\{v_3,u_4,u_5,u_6\}$ contains two edge-disjoint $3K_2$, say $M_1=\{w_2v_4,w_1u_1,u_2u_3\}$ and $M_2=\{w_2v_5,w_1u_2,u_1u_3\}$.
Assume that $c(v_3u_5)=c(e)$ for some $e\in M_1$. Then $\{v_3u_5,u_4u_6\}\cup M_2$ is a rainbow $5K_2$ in $T_n$. Similarly, $c(v_4u_2)=c(u_1u_3)$ and $c(v_5u_6)=c(u_4u_5)$. Since $e_G(S,A_i)\ge 5$ for each $i\in[2]$, we may assume $N_{A_1}(w_j)=V(A_1)$ and $N_{A_2}(w_k)=V(A_2)$, where $j,k\in[2]$. We have $j= k$, since otherwise, $\{w_ju_1,w_ku_4,v_3u_5,v_4u_2,v_5u_6\}$ is a rainbow $5K_2$ in $T_n$. Assume $j=k=1$. Then  $\{w_1u_4,w_2u_\ell,v_3u_5,v_4u_2,v_5u_6\}$ is a rainbow $5K_2$ in $T_n$, where $\ell\in\{1,3\}$, a contradiction.\medskip

Assume now $s=3$. Then $q=n-5$. Let $S=\{w_1, w_2, w_3\}$. Note that $|V_0|\ge 5$.
 Suppose $|B|\ne 0$. Then $|B|=2$ and $a_i=1$ for each $i\in[q]$. Hence, $e(G)=e(G[S])+e_G(S,V_0\cup B)+e(G[B])=3+(2n-4)+1<2n+1$, a contradiction. Thus $|B|=0$. Then $a_1=3$ and $a_i=1$ for each $i\ge 2$. Let $V(A_1)=\{u_1,u_2,u_3\}$. We claim that $e(G[S\cup V(A_1)])\le 3\cdot 6-7$. Suppose $e(G[S\cup V(A_1)])=3\cdot 6-6$. Then $e_G(S,V_0)=e(G)-e(G[S\cup V(A_1)])=2n-11$, which implies that $d_G(v_2)=3$ and so $G[S]=K_3$ since $G[S\cup V(A_1)]$ is a plane triangulation. Hence, $e_G(S,V(A_1))=6$ and $e_G(w_i,V(A_1))\ge 1$ for all $i\in[3]$ since $G[S]=K_3$. Thus, $d_G(v_i)=2$ for each $i\ge 3$ since otherwise $G$ contains a $K_{3,3}$-minor (with one part $\{w_1,w_2,w_3\}$ and the other part $\{v_2,v_3,V(A_1)\}$).
By Lemma~\ref{lem}(a), $G-\{v_i,v_j\}$ for $v_i,v_j\in V_0\less v_2$ contains two edge-disjoint $4K_2$ because $N_G(v_i)\cap N_G(v_j)\ne \emptyset$. By Claim, $E_{T_n}(V_0\less v_2)=\emptyset$.
Thus, $e_{T_n}(v_i,V(A_1)\cup\{v_2\})\ge 1$ because $\delta(T_n)\ge 3$. Since $|V_0|\ge 5$, there exist two vertices in $V_0\less v_2$ which has common neighbors in $G$. Assume that $N_G(v_3)=N_G(v_4)=\{w_2,w_3\}$. Then there exists one vertex in $\{w_2,w_3\}$, say $w_2$, such that $w_2v_5\in E(G)$ since $d_G(v_5)=2$.
Thus, $T_n$ contains a $K_{3,3}$-minor (with one part $\{w_2,\{w_1,w_3\},V(A_1)\cup\{v_2\}\}$ and the other part $\{v_3,v_4,v_5\}$), a contradiction. Thus, $e(G[S\cup V(A_1)])\le 3\cdot 6-7$. Then $2n+1=e(G)=e(G[S\cup V(A_1)])+e_G(S,V_0)\le (3\cdot 6-7)+(2(n-3)-4)=2n+1$, which implies that $e(G[S\cup V(A_1)])=11$ and $e_G(S,V_0)=2n-10$. This means $d_G(v_2)=d_G(v_3)=3$ and $d_G(v_i)=2$ for each $i\ge 4$. Thus, there are exactly two vertices in $S$, say $w_1,w_2$, having neighbors in $V(A_1)$, since otherwise $G$ contains a $K_{3,3}$-minor (with one part $\{w_1,w_2,w_3\}$ and the other part $\{v_2,v_3,V(A_1)\}$) or $e(G[S\cup V(A_1)])<11$. Note that $e_G(\{w_1,w_2\}, V(A_1))\ge5$. We next show $E_{T_n}(V_0)=\emptyset$. Suppose $E_{T_n}(V_0)\ne \emptyset$. Note that  $v_2v_3\notin E(T_n)$ because $w_1,w_2,w_3$ are in one face of $G[S\cup V(A_1)]$ and so $v_2$ and $v_3$ are in different faces of $G[S\cup V(A_1)]$, and $E_{T_n}(V_0\less\{v_2,v_3\})=\emptyset$ because $G-\{v_i,v_j\}$ for any $v_i,v_j\in V_0\less\{v_2,v_3\}$ contains two edge-disjoint $4K_2$. Hence, $E_{T_n}(\{v_2,v_3\},\{v_4,\ldots, v_q\})\ne \emptyset$. Assume $v_3v_4\in E(T_n)$. We show $c(v_3v_4)=c(v_2w_3)$. Suppose $c(v_3v_4)\ne c(v_2w_3)$. Then $G-\{v_3,v_4,v_2,w_3\}$ contains two edge-disjoint $3K_2$, say $M_1$ and $M_2$. Assume that $c(v_3v_4)=c(e)$ for some $e\in M_1$. Then $\{v_3v_4,v_2w_3\}\cup M_2$ is a rainbow $5K_2$ in $T_n$. Thus, $c(v_3v_4)=c(v_2w_3)$.
We next show $N_G(v_j)=\{w_1,w_2\}$ for $j\ge 5$. Suppose $w_3v_5\in E(G)$. Then $G-\{v_3,v_4\}$ contains two edge-disjoint $4K_2$. By Claim, $v_3v_4\notin E(T_n)$. Hence, $N_G(v_j)=\{w_1,w_2\}$ for $j\ge 5$. Observe that $e_{T_n}(v_i,V(A_1)\cup\{v_2,v_3\})\ge 1$ for each $i\ge 4$. Then $E_{T_n}(V(A_1)\cup\{v_2\},\{v_5,\ldots,v_q\})\ne \emptyset$ since otherwise $T_n$ contains a $K_{3,3}$-minor (with one part $\{v_4,v_5,v_6\}$ and the other part $\{\{w_1,w_3\},w_2,v_3\}$ or $\{w_1,\{w_2,w_3\},v_3\}$). Hence, either $v_jv_2\in E(T_n)$ or $v_ju_i\in E(T_n)$ for $i\in[3]$ and $j\ge 5$. Assume that $v_jv_2\in E(T_n)$.
If $c(v_jv_2)=c(v_3v_4)$, then $\{v_jv_2,v_3w_3,v_kw_2,w_1u_1,u_2u_3\}$ is a rainbow $5K_2$ in $T_n$, where $v_k\in V_0\less\{v_2,v_3,v_4,v_j\}$. Thus $c(v_jv_2)\ne c(v_3v_4)$. Then $G-\{v_2,v_3,v_4,v_j\}$ contains two edge-disjoint $3K_2$, say $M_1$ and $M_2$. Assume $c(v_jv_2)=c(e)$ for some $M_1$. Then $\{v_jv_2,v_3v_4\}\cup M_2$ is a rainbow $5K_2$, a contradiction. Hence, $v_ju_i\in E(T_n)$ for some $i\in[3]$ and $j\ge 5$. Assume $v_5u_1\in E(T_n)$. If $c(v_5u_1)=c(v_3v_4)$, then $\{v_5u_1,u_2u_3,v_3w_3,v_6w_2,w_1v_2,\}$ is a rainbow $5K_2$ in $T_n$. Thus $c(v_5v_2)\ne c(v_3v_4)$. Then $G-\{v_2,v_3,v_4,v_5\}$ contains two edge-disjoint $3K_2$, say $M_1$ and $M_2$. Assume $c(v_5v_2)=c(e)$ for some $M_1$. Then $\{v_5v_2,v_3v_4\}\cup M_2$ is a rainbow $5K_2$, a contradiction.
Hence, $E_{T_n}(V_0)=\emptyset$ and so $e_{T_n}(v_i,V(A_1))\ge 1$ for each $i\ge 4$.
If $v_jw_3\in E(G)$ for some $j\ge 4$, then $T_n$ contains a $K_{3,3}$-minor (with one part $\{w_1,w_2,\{w_3,v_j\}\}$ and the other part $\{v_2,v_3,V(A_1)\}$). Thus $N_G(v_i)=\{w_1,w_2\}$ for each $i\ge4$.
But then $T_n$ contains a $K_{3,3}$-minor (with one part $\{w_1,w_2,V(A_1)\}$ and the other part $\{v_4,v_5,v_6\}$) because $|V_0|\ge 5$, a contradiction.\medskip

Finally, we assume $s=4$. Let $S=\{w_1,w_2,w_3,w_4\}$. Then $q=n-4$ and $a_i=1$ for each $i\in[q]$. We claim $e(G[S])\le 5$. Suppose $e(G[S])=6$. Then $e_G(S,V_0)=e(G)-e(G[S])=2n-5$. We see $d_G(v_1)\le 3$ since otherwise $e(G[S\cup\{v_1\}])\ge 10$, which implies $d_G(v_2)=d_G(v_3)=3$ because $e_G(S,V_0)=2n-5$.
We next show $d_G(v_4)=3$. Suppose $d_G(v_4)=2$. Then $d_G(v_i)=2$ for each $i\ge 4$.
We next show $E_{T_n}(\{v_4,\ldots,v_q\})=\emptyset$. Suppose $v_4v_5\in E(T_n)$. Assume  $N_G(v_6)\cap N_G(v_7)\ne \emptyset$.  By Lemma~\ref{lem}(a), $G-\{v_4,v_5\}$ contains two edge-disjoint $4K_2$ and so $T_n$ has a rainbow $5K_2$. Thus $N_G(v_6)\cap N_G(v_7)=\emptyset$. But then $G-\{v_4,v_5\}$ contains two edge-disjoint $4K_2$ because $G[\{w_i,w_j,v_1,v_2,v_3\}]$ contains a $2K_2$ for any $1\le i<j\le4$, which implies $T_n$ has a rainbow $5K_2$.
Since $G[S]$ has four faces, there is at most one vertex in $\{v_4,\ldots,v_q\}$, say $v_q$, such that  $e_{T_n}(v_q,S)=3$.
Thus $e_{T_n}(v_i,\{v_1,v_2,v_3\})\ge1$ for each $4 \le i \le q-1$ because $\delta(T_n)\ge 3$.
Without loss of generality, assume that $ N_G(v_1) =\{w_1,w_2,w_3\}$, $N_G(v_2)=\{w_1,w_2,w_4\}$ and $N_G(v_3)=\{w_1,w_3,w_4\}$. If $v_4,v_5,v_6$ are adjacent to the same vertex of $\{v_1,v_2,v_3\}$, say $v_1$, then $T_n$ has a rainbow $5K_2$ since $G-\{v_4,v_1\}$ contains two edge-disjoint $4K_2$, a contradiction. Hence, we can find two vertices of $\{v_4,v_5,v_6\}$, say $v_4$ and $v_5$, which are adjacent to different vertices of $\{v_1,v_2,v_3\}$. Without loss of generality, assume that $v_1v_4\in E(T_n)$ and $v_2v_5 \in E(T_n)$. Let $c(v_1v_4)=c(e_1)$ and $c(v_2v_5)=c(e_2)$ for $e_i\in E(G)$. We claim $e_1\ne e_2$. Suppose $e_1=e_2$. If $e_1=e_2\in G[S]$, then $\{v_1v_4,v_2w_i,v_3w_j,v_5w_k,v_6w_\ell\}$ is a rainbow $5K_2$, where $\{i,j,k,\ell\}=[4]$. Hence, $e_1=e_2=e\notin G[S]$. Assume $e=v_3w_4$. Then $G-\{v_1,v_4,v_3,w_3\}$ has a $3K_2=M_1$. This implies  $M_1\cup\{v_1v_4,v_3w_3\}$ is a rainbow $5K_2$, a contradiction. Assume $e\ne v_3w_4$. Then $G-\{v_1,v_4,v_3,w_4\}-e$ has a $3K_2$ when $e=v_1w_i$ or $e=v_jw_k$ for $j\ge 5$ and some $i,k\in[4]$ and $G-\{v_2,v_5,v_3,w_4\}-e$ has a $3K_2$ when $e=v_2w_i$ or $e=v_jw_k$ for $3\le j\le 4$ and some $i,k\in[4]$. This implies $T_n$ has a rainbow $5K_2$, a contradiction. Hence, $e_1\ne e_2$. Moreover, we may further assume each matching of size $2$ in $T_n[V_0]$ has distinct colors. If $e_1,e_2\in G[S]$ or $e_1,e_2\notin G[S]$, then $\{v_1v_4,v_2v_5,v_3w_i,v_6w_j,v_7w_k\}$ is a rainbow $5K_2$ in $T_n$, where $i,j,k\in[4]$. Hence, only one of $\{e_1,e_2\}$ belongs to $G[S]$. Without loss of generality, assume $e_1\in G[S]$ and $e_2\notin G[S]$. Then $e_2=v_3w_i$ and $N_G(v_6)=N_G(v_7)=N_G(v_3)\less w_i$ since otherwise $T_n$ has a rainbow $5K_2$, which implies either $v_iv_3\in E(T_n)$ for some $i\in\{6,7\}$. If $c(v_iv_3)=c(e)$ for some $e\in G[S]$, then $e_1,e\in G[S]$ and so $T_n$ has a rainbow $5K_2$ by using the above method. Hence, $c(v_iv_3)=c(e')$ for some $e'\notin G[S]$. But then $e_2,e'\notin G[S]$ and so $T_n$ has a rainbow $5K_2$, a contradiction.
Thus $d_G(v_4)=3$ and so $d_G(v_5)=\cdots=d_G(v_{q-1})=2$ and $d_G(v_q)=1$. By Claim, $E_{T_n}(\{v_5,\ldots,v_q\})=\emptyset$ since $G[S\cup \{v_1,\ldots,v_4\}]$ contains two edge-disjoint $4K_2$. Thus, $E_{T_n}(\{v_1,\ldots,v_4\},v_i)\ne \emptyset$ for each $i\ge 5$ because $\delta(T_n)\ge 3$.
Assume $v_4v_5\in E(T_n)$. Let $N_G(v_6)=\{w_3,w_4\}$. By Lemma~\ref{lem}(a), $G-\{v_4,v_5,v_6,w_i\}$ for each $i\in\{3,4\}$ contains two edge-disjoint $3K_2$. Hence, $c(v_4v_5)=c(v_6w_i)$ for each $i\in\{3,4\}$ since otherwise $T_n$ has a rainbow $5K_2$. But then $c(v_6w_3)\ne c(v_6w_4)$, a contradiction.
Thus, $e(G[S])\le 5$.
Then $2n+1=e(G)=e(G[S])+e_G(S,V_0)\le 5+2n-4=2n+1$, which implies that $e_G(S,V_0)=2n-4$ and $e(G[S])=5$. We claim that $d_G(v_1)\le 3$. Suppose  $d_G(v_1)=4$. Then $d_G(v_2)=d_G(v_3)=3$ and $d_G(v_i)=2$ for all $i\ge4$.
Note that $G-\{v_i,v_j\}$ contains two edge-disjoint $4K_2$ for $4\le i<j\le q$. By Claim, $E_{T_n}(\{v_4,\ldots,v_q\})=\emptyset$. Hence, $e_{T_n}(v_i,\{v_1,v_2,v_3\})\ge1$ for each $i\ge 4$.
Assume $e_{T_n}(v_1,v_i)\ge 1$ for each $i\in\{4,5,6,7\}$. Then each vertex of $S$ has neighbors in $\{v_4,\ldots,v_q\}\less v_i$ for $i\in
\{4,5,6,7\}$. Hence, $G-\{v_1v_i\}$ has two edge-disjoint $4K_2$ because $G[\{v_2,v_3,w_i,w_j\}]$ has a $2K_2$ for $1\le i<j\le 4$. By Claim, $v_1v_i\notin E(T_n)$, a contradiction. Thus, there exists two vertices of $\{v_4,v_5,v_6,v_7\}$, say $v_4$ and $v_5$, having distinct neighbors in $\{v_1,v_2,v_3\}$. Using the above similar method, one can check $T_n$ has a rainbow $5K_2$, a contradiction.
Thus $d_G(v_1)=\cdots=d_G(v_4)=3$ and so $d_G(v_i)=2$ for each $i\ge5$. Note that $G[S]$ has only three faces.
Then either $N_G(v_{i'})\ne N_G(v_{j'})$ for any $1\le i'<j'\le 4$ or $N_G(v_i)=N_G(v_j)$ and $N_G(v_k)=N_G(v_\ell)$ for some four distinct integers $i,j,k,\ell\in[4]$.
It is easily to check that either $G[S\cup\{v_1,\ldots,v_4\}]$ or $G-\{v_i,v_j\}$ contains two edge-disjoint $4K_2$, where $i\in[4]$ and $j\ge 5$. By Claim, $E_{T_n}(\{v_5,\ldots,v_q\})=\emptyset$ and $E_{T_n}(\{v_5,\ldots,v_q\},\{v_1,\ldots,v_4\})=\emptyset$. Hence, $e(T_n)=e(G[S\cup\{v_1,v_2,v_3,v_4\}])+e_G(S,\{v_5,\ldots,v_q\})\le 3\cdot 8-6+2(n-8)=2n+2<3n-6$, a contradiction.

The proof is thus complete.\qed\\

\noindent{\bf Acknowledgments.}
The authors would like to thank the anonymous referees and the Associate Editor Florian Pfender for many valuable suggestions and comments. Z. Qin was partially supported by the Fundamental Research Funds for the Central Universities (No. 300102128104). Y. Lan and Y. Shi were partially supported by the National Natural Science
Foundation of China (No.  11771221 and 11811540390), Natural Science Foundation of Tianjin (No. 17JCQNJC00300)  and  China--Slovenia bilateral
project ``Some topics in modern graph theory" (No. 12-6).

\end{document}